\documentclass[9pt]{article} 
\date{}
\usepackage{amssymb}
\usepackage{amsmath}
\usepackage{amsthm}
\usepackage{cite}
\usepackage{abstract}
\usepackage{graphicx}
\allowdisplaybreaks[4]
\begin{document}
	\title{The Existence of Pseudoharmonic Maps For Small Horizontal Energy \footnote{Supported by NSFC Grants No.11771087 and No.12171091.} } 
		\author{Biqiang Zhao}
	\maketitle

	\setlength{\abstitleskip}{-2em}
	\renewcommand{\abstractname}{}		
	\begin{abstract}
	$\mathbf{Abstract:}$ In this paper, we consider the pseudoharmonic heat flow with small initial horizontal energy and give the existence of pseudoharmonic maps from closed pseudo-Hermitian manifolds to closed Riemannian manifolds.
	\end{abstract}
    \section{Introduction} 
    \par   
    \ \ \ The existence of harmonic mappings is an important issue which has been widely studied in Geometric Analysis. Eells and Sampson \cite{MR164306} proved the existence of harmonic maps from closed Riemannian manifolds into nonpositively curved closed Riemannian manifolds through the heat flow method. As for general target manifolds, we cannot get an Eell-Sampson's type theorem in general. In \cite{MR345129}, Mitteau discussed the harmonic heat flow with small initial energy and established an existence theorem of harmonic maps when target manifold $N$ is not assumed to be nonpositively curved. \par 
    For pseuduoharmonic case, Barletta $et\ al.$ in \cite{MR1871387} introduced the notion of pseudoharmonic maps in CR geometry, which is a generalization of harmonic maps.  Chang and Chang \cite{MR3038722} studied the solution of the pseudoharmonic heat flow and proved the existence of pseudoharmonic maps from closed pseudo-Hermitian manifolds into a nonpositively curved closed Riemannian manidfolds under the commutation condition $[\Delta_b,\xi]=0$, where $\Delta_b$ is the sub-Laplacian and $\xi$ is the Reeb vector field. In \cite{MR3844509}, Ren and Yang obtained the Eells-Sampson's type theorem without the commutation condition in \cite{MR3038722}.\par 
    In the present paper, following the idea of Mitteau, we obtain the global existence of the pseudoharmonic heat flow 
    and the existence of pseudoharmonic maps if the initial horizontal energy is small enough. The main result is the following theorem.

    ~\\
    $\mathbf{Theorem\ 1.1}$ Let ($M^{2m+1},HM,J,\theta$) be a closed pseudo-Hermitian manifold with $m\geq 2$. Let $(N^n,g_N)$ be a closed Riemannian manifold and $\kappa$ be an upper bound for the sectional curvature of $N$. Suppose $h:M\rightarrow N$ is a smooth map with $e(h)\leq D$ for some constant $D\geq$0. Then there exists a constant $\epsilon\geq 0$ depending on $M,N,D$ such that if the horizontal energy $E_b(h)\leq \epsilon$, then the pseudoharmonic heat flow
    \begin{eqnarray}
    \frac{\partial u}{\partial t}=\tau(u),u(x,0)=h(x) 
    \end{eqnarray}
    has a solution on $M\times[0,\infty)$. Moreover, there exists a sequence $t_i\rightarrow \infty$ such that $u(\cdot,t_i)\rightarrow u_\infty$ uniformly to a pseudoharmonic map $u_{\infty}$.
    
    ~\par
    Assume that $u$ is a solution of the pseudoharmonic heat flow (1), then the key point now is to estimate the supremum of the total energy density of $u$, which will be denoted by $\rho$ in the following passage. The proof is mainly divided into two parts. First we point out that $\rho$  has an upper bound in a finite time interval $[0,t_0]$ through a control function. Next  for $t\geq t_0$, we show that if the initial horizontal energy is sufficiently small, $\rho$ has an upper bound by the continuity of $\rho$ of t.    
    \\

    \section{Pseudo-Hermitian Geometry}
    In this section we introduce some basic notations in pseudo-Hermitian geometry (cf. \cite{MR2214654,MR520599,MR0399517}). A smooth  manifold $M^{2m+1}$ is a CR manifold if there is a smooth rank $m$ complex subbundle $T^{1,0}M\subseteq TM\otimes \mathbb{C}$  
    such that  
    \begin{eqnarray}
    &&T^{1,0}M\cap T^{0,1}M=\{0\} \nonumber\\
    &&[T^{1,0}M,T^{1,0}M]\subseteq T^{1,0}M \nonumber
    \end{eqnarray}
    where $T^{0,1}M=\overline{T^{1,0}M}$. Equivalently, the CR structure may also be described by the real bundle $HM=Re\{T^{1,0}M\oplus T^{0,1}M\}$ and an almost complex structure $J$ on $HM$, $J(X+\bar{X})=\sqrt{-1}(X-\bar{X}), \forall X\in T^{1,0}M$. Set
    \begin{eqnarray}
    E=\{\omega\in T^*M|\omega(HM)=0\}. \nonumber
    \end{eqnarray}
    If $M$ is oriented, then $E$ admits global nowhere vanishing sections. A section $\theta \in\Gamma(E \backslash \{0\})$ is called a pseudo-Hermitian structure. The Levi form is defined by
    \begin{eqnarray}
    L_\theta(X,Y)=d\theta(X,JY),\ \forall X,Y\in HM. \nonumber
    \end{eqnarray}
    A pseudo-Hermitian manifold is the quadruple $(M,HM,J,\theta)$ with $L_\theta$ is positive definite.\par
    The Reeb vector filed is the unique vector field $\xi$ on M such that 
    \begin{eqnarray}
    \theta(\xi)=1,\ d\theta(\xi,\cdot)=0. \nonumber
    \end{eqnarray}
    Then there is a decomposition of the tangent bundle $TM$, $TM=HM\oplus \mathbb{R}\xi$, which induces the projection $\pi_b:TM\rightarrow HM$. We can extend $L_\theta$ to get a Riemannian metric $g_\theta$, called Webster metric, by
    \begin{eqnarray}
    g_\theta=L_\theta+\theta\otimes \theta. \nonumber
    \end{eqnarray}
    On a pseudo-Hermitian manifold, there exists a canonical connection $\nabla$, called the Tanaka-Webster connection (cf.\cite{MR520599}), such that 
    \begin{eqnarray}
    &&(1)\ \nabla_X\Gamma(HM)\subseteq \Gamma(HM),\ \forall X\in TM; \nonumber\\
    &&(2)\ \nabla g_\theta=0,\ \nabla J=0; \nonumber\\
    &&(3)\ T_\nabla(X,Y)=2d\theta(X,Y)\xi,\ T_\nabla(\xi,JX)+JT_\nabla(\xi,X)=0,\ \forall X,Y\in HM. \nonumber
    \end{eqnarray}
    The pseudo-Hermtian torsion $A$ is defined by $A(X,Y)=g_\theta(T_\nabla(\xi,X),Y)$ for any $X,Y\in TM$.\par 
    Let $(M,HM,J,\theta)$ be a pseudo-Hermitian manifold of dimension $2m+1$. Assume that $\{\eta_1,\cdots,\eta_m\}$ be a local frame of $T^{1,0}M$ on an open set $U\subseteq M$ and $\{\theta^1,\cdots,\theta^m\}$ its dual coframe. Then we have the structure equations for the Tanaka-Webester connection.
    \begin{eqnarray}
    &&d\theta=2\sqrt{-1}\theta^\alpha \wedge \theta^{\bar{\alpha}}, \nonumber\\
    &&d\theta^{\alpha}=\theta^\beta\wedge \theta^\alpha_\beta+A_{\bar{\alpha}\bar{\beta}}\theta\wedge\theta^\beta,\ \theta^\alpha_\beta+\theta^{\bar{\beta}}_{\bar{\alpha}} =0,\nonumber\\
    &&d\theta^\alpha_\beta=\theta^\gamma_\beta\wedge\theta^\alpha_\gamma+\Pi^\alpha_\beta \nonumber
    \end{eqnarray}
    with
    \begin{eqnarray}
    \Pi^\alpha_\beta=2\sqrt{-1}(\theta^\alpha\wedge\tau_\beta+\theta_\beta\wedge\tau^\alpha)+R^\alpha_{\beta\lambda\bar{\mu}}\theta^\lambda\wedge\theta^{\bar{\mu}}+W^\alpha_{\beta\bar{\mu}}\theta\wedge\theta^{\bar{\mu}}-W^\alpha_{\beta\mu}\theta\wedge\theta^{\bar{\mu}} \nonumber
    \end{eqnarray}   
    where $\{\theta_\beta^\alpha\}$ are the Tanaka-Webster connection 1-forms with respect to $\{\eta_\alpha\}$, $W^\alpha_{\beta\bar{\mu}}=A^\alpha_{\bar{\mu}, \beta},\ W^\alpha_{\beta\mu}=A^{\bar{\mu}}_{\beta,\bar{\alpha}},\ R^\alpha_{\beta\lambda\bar{\mu}}$ is the Webster curvature. 
    
    \par
    For a Riemmanian manifold $(N,g_N)$ with the Levi-Civita connection $\nabla^N$. Let $\{\sigma^i\}$ be a local orthonnormal coframe of $T^*N$ and $\{E_i\}$ its dual frame of $TN$. For a smooth map $f:M\rightarrow N$, there are two induced connections on $f^*TN$ and $TM\otimes f^*TN$, also denoted by $\nabla$. Now we give some basic notations of energy and energy density. Let us recall the total energy density $e(f)$ and the total energy $E(f)$ for a $C^2$ map $f:M\rightarrow N$ 
    \begin{eqnarray}
    e(f)=\frac{1}{2}|df|^2,\ E(f)=\frac{1}{2}\int_Me(f) dV \nonumber
    \end{eqnarray}
    where $dV=\theta \wedge (d\theta)^m$ is the volume form. For pseudo-Hermitian geometry, we can define the horizontal energy  density $e_b(f)$ and the horizontal energy $E_b(f)$ by
    \begin{eqnarray}
    e_b(f)=\frac{1}{2}|d_bf|^2=\frac{1}{2}|df\circ\pi_b|^2,\ E_b(f)=\frac{1}{2}\int_Me_b(f)dV .\nonumber
    \end{eqnarray}
    We also consider the vertical energy density $e_0(f)$ and the vertical energy $E_0(f)$
    \begin{eqnarray}
    e_0(f)=e(f)-e_b(f)=\frac{1}{2}|f_0|^2,\ E_0(f)=\frac{1}{2}\int_Me_0(f)dV.\nonumber
    \end{eqnarray}
    
    $\mathbf{Definition\ 2.1}$ (cf.\cite{MR1871387}) A smooth map $f:M\rightarrow N$ is called pseudoharmonic if the tensor field
    \begin{eqnarray}
    \tau(f)=0, \nonumber
    \end{eqnarray}
    where $\tau(f)=trace_{g_\theta}(\nabla_bd_bf),\  \nabla_bd_bf=\nabla d_bf\circ\pi_b,\ d_bf=df\circ\pi_b$.
    \par  
    It is proved in \cite{MR1871387} that the pseudoharmonic map is the critical point of the horizontal energy. In order to prove the existence of pseudoharmonic maps, as in Riemannian case, we consider the pseudoharmonic heat flow
    \begin{eqnarray}
    \frac{\partial u}{\partial t}=\tau(u),u(x,0)=h(x).  \nonumber
    \end{eqnarray}  
    For short time existence, one can refer to \cite{MR3844509}. Therefore we are only concerned in this paper with the long time existence of the pseudoharmonic heat flow with small initial horizontal energy. Now we give some results which will be used in the proof of Theorem 1.1. In this paper, we employ the index conventions 
    \begin{eqnarray}
    A,B,C&=&0,1,\cdots,m,\bar{1},\cdots,\bar{m}, \nonumber\\
    \alpha,\beta,\gamma&=&1,2,\cdots,m, \nonumber\\
    i,j,k&=&1,2,\cdots,n, \nonumber
    \end{eqnarray}
    and use the Einstein summation convention. First let us recall the CR Bochner formulas.
    
    ~\\
    $\mathbf{Lemma\ 2.1 }$ (cf.\cite{MR3194214})  We chooose a local coframe $\{\theta,\theta^\alpha,\theta^{\bar{\alpha}}\}$ on $M$, a local frame $\{E_j\}$ on $N$ and denote the pseudo-Hermitian Ricci curvature by $R^M_{\alpha\bar{\beta}}$. For any smooth map $f:M\rightarrow N$, we denote the components of $df$, $\nabla df$ by $f^i_A,\ f^i_{AB}$. Then we have 
    \begin{eqnarray}
    \Delta_b|d_bf|^2&=&2|\nabla_bd_bf|^2+2\langle \nabla_b\tau (f),d_bf \rangle +8\sqrt{-1}(f^i_{\bar{\alpha}}f^i_{0\alpha}-f^i_{\alpha}f^i_{0\bar{\alpha}})  \nonumber\\
    &&+4R^M_{\alpha\bar{\beta}}f^i_{\bar{\alpha}}f^i_\beta-4\sqrt{-1}(m-2)(f^i_\alpha f^i_\beta A_{\bar{\alpha}\bar{\beta}}-f^i_{\bar{\alpha}}f^i_{\bar{\beta}}A_{\alpha\beta}) \nonumber\\
    &&+4(f^i_{\bar{\alpha}}f^j_{\beta}f^k_{\bar{\beta}}f^i_{\alpha}R^N_{ijkl}+f^i_\alpha f^j_\beta f^k_{\bar{\alpha}}f^i_{\bar{\beta}}R^N_{ijkl})\\
    \Delta_b|f_0|^2&=&2|\nabla_bf_0|^2+2\langle \nabla_\xi\tau(f),f_0\rangle+4f^i_0f^j_\alpha f^k_{\bar{\alpha}}f^l_0R^N_{ijkl} \nonumber\\
    &&+4(f^i_0f^i_\beta A_{\bar{\beta}\bar{\alpha},\alpha}+f^i_0f^i_{\bar{\beta}}A_{\beta\alpha,\bar{\alpha}}+f^i_0f^i_{\beta\alpha}A_{\bar{\beta}\bar{\alpha}}+f^i_0f^i_{\bar{\beta}\bar{\alpha}}A_{\beta\alpha}).
    \end{eqnarray}  
    
    ~\\
    $\mathbf{Corollary\ 2.1}$ (cf.\cite{MR3844509}) There exists contansts $C_1,C_2$, where $C_1$ only depending on the pseudo-Hermitian Ricci curvature, the pseudo-Hermitian torsion and its divergence and $C_2$ only depending on the sectional curvature of $N$, such that for a solution $u$ of (1), we have 
    \begin{eqnarray}
    (\Delta_b-\partial_t)e(u)\geqslant -C_1e(u)-C_2e^2(u).
    \end{eqnarray}
    In order to estimate the vertical energy, we needs some formulas in \cite{MR3938843}.

    ~\\
    $\mathbf{Lemma\ 2.2}$ (cf.\cite{MR3938843}) Let ($M^{2m+1},HM,J,\theta$) be a closed pseudo-Hermitian manifold with $m\geq 2$ and $N^n$ be a closed Riemannian manifold. Suppose $f:M\rightarrow N$ is a smooth map. Then
    \begin{eqnarray}
    \sqrt{-1}(f^i_\alpha f^i_{\bar{\alpha}0}-f^i_{\bar{\alpha}}f^i_{\alpha0})&=&\frac{1}{m}\langle Pf+\overline{Pf},d_bf \rangle-\frac{1}{2m}\langle d_bf,\nabla_b\tau(f)\rangle \nonumber\\
    &&+\sqrt{-1}(f^i_\alpha f^i_\beta A_{\bar{\alpha}\bar{\beta}}-f^i_{\bar{\alpha}}f^i_{\bar{\beta}}A_{\alpha\beta})  \\
    \int_M \langle Pf+\overline{Pf},d_bf \rangle dV &\leq& \frac{2m}{m-1}\int_Mf^i_\alpha f^j_\beta f^k_{\bar{\alpha}}f^l_{\bar{\beta}}R^N_{jikl}dV 
    \end{eqnarray}
    where $Pf=(P^j_\beta f)\theta^\beta\otimes E_j,\ P^j_\beta f=f^j_{\bar{\alpha}\alpha\beta}+2\sqrt{-1}mA_{\beta\alpha}f^j_{\bar{\alpha}}$.
    \\
    
    $Remark\ 2.2$ (6) is contained in the proof of Theorem 4.1 of \cite{MR3938843}.
    ~\\
    
    $\mathbf{Lemma\ 2.3}$ (cf. \cite{MR3844509}) If $\phi\in C^\infty(M\times (0,\delta))$ is nonnegative and satisfies
    \begin{eqnarray}
    (\Delta_b-\partial_t)\phi\ge 0, \nonumber
    \end{eqnarray}
    then for any $\epsilon\in(0,\delta),t\in[\epsilon,\delta)$, we find that
    \begin{eqnarray}
    \phi(x,t)\le C_\epsilon\int_{t-\epsilon}^{t}\int_M\phi(y,s)dVds, \nonumber
    \end{eqnarray}
    where $C_\epsilon$ only depends on $\epsilon$.

    \section{Long Time Existence for Small Horizontal Energy}
   Let $u=u(x,t)$ be a solution of (1) and $[0,T_0)$ be the maximal existence time interval of $u$, where $0\le T_0\leq +\infty$. Then we have the following lemma.

   ~\\
   $\mathbf{Lemma\ 3.1}$ If $ e(h)\leq D$ in the pseudoharmonic heat flow (1), then for the maximal time $T_0$ of (1) and the total energy density $e(u)$, we have the following estimate 
   \begin{eqnarray}
   T_0 \geq\frac{1}{C_1}log(1+\frac{C_1}{DC_2}) 
   \end{eqnarray}
   and
   \begin{eqnarray}
   e(u)\leq \frac{C_1De^{C_1t}}{C_1+C_2D-C_2De^{C_1t}} 
   \end{eqnarray}
   on $[0,\frac{1}{C_1}log(1+\frac{C_1}{DC_2})\ ) $, where $C_1,C_2$ are given by Corollary 2.1.
   
   ~\\
   $\mathbf{Proof. }$ From Corollary 2.1, we have 
   \begin{eqnarray}
   (\Delta_b-\partial_t)e(u)\geqslant -C_1e(u)-C_2e^2(u). \nonumber
   \end{eqnarray}
   We define a function $g$ by 
   \begin{eqnarray}
   g(t)&=&\frac{C_1}{C_2}(\frac{1}{1-C_2e^{B+C_1t}}-1)=\frac{C_1}{C_2}\frac{C_2e^{B+C_1t}}{1-C_2e^{B+C_1t}}\nonumber\\
   &=&\frac{C_1e^{B+C_1t}}{1-C_2e^{B+C_1t}} \nonumber
   \end{eqnarray}
   and $g$ satisfies 
   \begin{eqnarray}
   \frac{\partial g}{\partial t}&=&\frac{C_1}{C_2}\frac{C_2C_1e^{B+C_1t}}{(1-C_2e^{B+C_1t})^2}=\frac{C_1^2e^{B+C_1t}}{(1-C_2e^{B+C_1t})^2} \nonumber\\
   &=&g\cdot \frac{C_1}{1-C_2e^{B+C_1t}}=g(C_1+C_2g) >0,
   \end{eqnarray}
   where  $B=log\frac{D+\delta}{C_1+C_2(D+\delta)},\ \delta>0$. Also 
   \begin{eqnarray}
   g(0)&=&\frac{C_1e^B}{1-C_2e^{B}}=\frac{C_1\frac{D+\delta}{C_1+C_2(D+\delta)}}{1-C_2\frac{D+\delta}{C_1+C_2(D+\delta)}} \nonumber\\
   &=&\frac{C_1(D+\delta)}{C_1}=D+\delta >e(h).\nonumber
   \end{eqnarray}
   Now we claim that:
   \begin{eqnarray}
   e(u)(x,t) < g(t),\ \forall(x,t)\in M\times[0,\frac{-(B+logC_2)}{C_1}) .
   \end{eqnarray}
   Otherwise, there exits a $T\in (0,\frac{-(B+logC_2)}{C_1})$ is the first time such that $\inf\limits_{M\times[0,T]}(g-e(u))=0$. Note that 
   \begin{eqnarray}
   (\triangle_b-\partial_t)(g-e(u))&\leq& -C_1g-C_2g^2+C_1e(u)+C_2e(u) \nonumber\\
   &=&-C_1(g-e(u))-C_2(g-e(u))(g+e(u)) .
   \end{eqnarray}
   From (9), we have 
   \begin{eqnarray}
   g(t)+e(u)>g(0)=D+\delta  
   \end{eqnarray}
   on $(0,T)$. (11)+(12) yield that
   \begin{eqnarray}
   (\triangle_b-\partial_t)(g-e(u))\leq -(C_1+C_2(D+\delta))(g-e(u)) \nonumber
   \end{eqnarray}
   and
   \begin{eqnarray}
   (\triangle_b-\partial_t)e^{-(C_1+C_2(D+\delta))t}(g-e(u))\leq 0. \nonumber
   \end{eqnarray}
   By the maximum principle, we have 
   \begin{eqnarray}
   e^{-(C_1+C_2(D+\delta))t}(g-e(u))\geq g(0)-e(h)\geq \delta, \nonumber
   \end{eqnarray}
   $\forall x\in M,\ t\in [0,T]$. Let $t=T$, we find that
   \begin{eqnarray}
   (g-e(u))(x,T)\geq \delta e^{(C_1+C_2(D+\delta))T} >0,\ \forall x\in M .\nonumber
   \end{eqnarray}
   But this leads a contradiction with $\inf\limits_{M\times[0,T]}(g-e(u))=0$, which proves the claim. We can complete the proof by taking $\delta\rightarrow 0$ in (3.4). \\
   {\qed}
   
   ~\\ 
   $\mathbf{ Proof\ of\ Theorem\ 1.1}$ By a similar argument of the short-time existence in \cite{MR3844509}, the long-time existence for (1) will hold if $e(u)$ is uniformly bounded inside any finite time. Due to Lemma 3.1 and (9), there exists a $t_0$ such that 
   \begin{eqnarray}
   e(u)(x,t)\leq g(t)\leq g(t_0)= 2D 
   \end{eqnarray}
   on [0,$t_0$].\par 
   For a given $T\geq t_0$, we denote
   \begin{eqnarray}
   \rho=\sup\limits_{M\times[0,T]}e(u)(x,t) .\nonumber
   \end{eqnarray}
   Then Corollary 2.1 yields that 
   \begin{eqnarray}
   (\triangle_b-\partial_t)e(u)\geq -(C_1+C_2\rho)e(u), \nonumber
   \end{eqnarray}
   which implies 
   \begin{eqnarray}
   (\triangle_b-\partial_t)e^{-(C_1+C_2\rho)t}e(u) \geq 0 \nonumber
   \end{eqnarray}
   on $M\times[0,T]$. Then Lemma 2.3 yields that 
   \begin{eqnarray}
   &&e^{-(C_1+C_2\rho)t}e(u)(x,t)\nonumber\\
   &\leq& C_s\int_{t-s}^{t}\int_Me^{-(C_1+C_2\rho)r}e(u)(y,r)dV_ydr \nonumber\\
   &\leq&C_se^{-(C_1+C_2\rho)(t-s)}	\int_{t-s}^{t}\int_M[e_b(u)+e_0(u)](y,r)dV_ydr,
   \end{eqnarray}
   $\forall x\in M,\ 0<s\leq t\leq T$, where $C_s$ only depends on $s$. Integtate (5), we conclude that   
   \begin{eqnarray}
   &&\int_{t-s}^{t}\int_M\sqrt{-1}(u^i_\alpha u^i_{0\bar{\alpha}}-u^i_{\bar{\alpha}}u^i_{0\alpha})(y,r)dV_ydr\nonumber\\   
   &=&\frac{1}{m}\int_{t-s}^{t}\int_M\langle Pu+\overline{Pu},d_bu \rangle (y,r)dV_ydr\nonumber \\
   &&+2\int_{t-s}^{t}\int_M\sqrt{-1}(u^i_\alpha u^i_\beta A_{\bar{\alpha}\bar{\beta}}-u^i_{\bar{\alpha}}u^i_{\bar{\beta}}A_{\alpha\beta})  (y,r)dV_ydr\nonumber\\
   && -\frac{1}{2m}\int_{t-s}^{t}\int_M\langle d_bu,\nabla_b\tau(u) \rangle (y,r)dV_ydr.
   \end{eqnarray}
   The left hand side of (15) is 
   \begin{eqnarray}
   &&\int_{t-s}^{t}\int_M\sqrt{-1}(u^i_\alpha u^i_{0\bar{\alpha}}-u^i_{\bar{\alpha}}u^i_{0\alpha})
   \nonumber\\
   &=&\int_{t-s}^{t}\int_M\sqrt{-1}u^i_0(u^i_{\bar{\alpha}\alpha}-u^i_{\alpha\bar{\alpha}})(y,r)dV_ydr \nonumber\\
   &=&2m\int_{t-s}^{t}\int_M |u^i_0|^2 (y,r)dV_ydr .
   \end{eqnarray}
   Using (6), the terms in the right hand side of (15) can be estimated as
   \begin{eqnarray}
   &&\frac{1}{m}\int_{t-s}^{t}\int_M\langle Pu+\overline{Pu},d_bu \rangle(y,r)dV_ydr \nonumber\\
   &\leq& \frac{2}{m-1}\kappa \rho\int_{t-s}^{t}\int_M |d_bu|^2(y,r)dV_ydr ,
   \end{eqnarray}
   \begin{eqnarray}
   &&-\frac{1}{2m}\int_{t-s}^{t}\int_M\langle d_bu,\nabla_b\tau(u)\rangle (y,r)dV_ydr  \nonumber \\
   &=& -\frac{1}{4m}\int_{t-s}^{t}\int_M \partial_r \langle d_bu,d_bu \rangle (y,r)dV_ydr \nonumber\\
   &\leq& \frac{1}{2m}(E_b(u)(t-s)-E_b(u)(t)),
   \end{eqnarray}
   \begin{eqnarray}
   &&2\int_{t-s}^{t}\int_M\sqrt{-1}(u^i_\alpha u^i_\beta A_{\bar{\alpha}\bar{\beta}}-u^i_{\bar{\alpha}}u^i_{\bar{\beta}}A_{\alpha\beta}) (y,r)dV_ydr\nonumber\\
   &\leq& 2\kappa^\prime\int_{t-s}^{t}\int_M |d_bu|^2(y,r)dV_ydr,
   \end{eqnarray}   
   where $K^N\leq \kappa,\ |A|\leq \kappa^\prime$. Since 
   \begin{eqnarray}
   \frac{d}{dr}E_b(u)(r)&=&\int_M\langle \nabla_{\partial_r}d_bu,d_bu\rangle(y,r)dV_y=\int_M\langle \nabla_{b}\partial_ru,d_bu\rangle(y,r)dV_y \nonumber\\
   &=&-\int_M|\partial_ru|^2(y,r)dV_y \leq 0,
   \end{eqnarray}
   we have
   \begin{eqnarray}
   E_b(u)\leq E_b(h)\leq \epsilon. 
   \end{eqnarray}
   Then (15)-(21) yields that 
   \begin{eqnarray}
   \int_{t-s}^{t}\int_Me_0(u)(y,r)dV_ydr\leq C_3(1+\rho)E_b(h)\leq C_3(1+\rho)\epsilon ,
   \end{eqnarray}
   where $C_3$ is a constant depends on $M,N,s$.
   Hence (14)-(22) asserts that 
   \begin{eqnarray}
   e^{-(C_1+C_2\rho)t}e(u)(x,t)\leq C_4(1+\rho)e^{-(C_1+C_2\rho)(t-s)}\epsilon, \nonumber 
   \end{eqnarray} 
   implying
   \begin{eqnarray}
   e(u)(x,t)\leq C_4(1+\rho)e^{(C_1+C_2\rho)s}\epsilon,
   \end{eqnarray}
   $\forall x\in M,\ 0<s\leq t\leq T$, where $C_4$ is a constant depends on $M,N,s$. We now are ready to estimate the total energy density $e(u)$. Let us choose $s$ sufficiently small such that 
   \begin{eqnarray}
   s<\frac{1}{D(4D+2)C_2}.
   \end{eqnarray}
   We now claim that 
   \begin{eqnarray}
   \rho< -\frac{1}{2}+\sqrt{\frac{1}{4}+\frac{1}{C_2s}}, \nonumber
   \end{eqnarray}
   if $\epsilon$ is small enough. We see this by examining two cases.
   \\
   (i) If $\rho=\sup\limits_Me(u)(x,t)$ for some $t\in[0,t_0]$, then (13) and (24) deduce that 
   \begin{eqnarray}
   \frac{1}{4}+\frac{1}{C_2s} >\frac{1}{4}+(4D+2)D >(2D+\frac{1}{2})^2 ,   \nonumber
   \end{eqnarray} 
   implying
   \begin{eqnarray}
   \rho\leq 2D< -\frac{1}{2}+\sqrt{\frac{1}{4}+\frac{1}{C_2s}}. \nonumber
   \end{eqnarray} 
   (ii) If $\rho=\sup\limits_{M\times [t_0,T]}e(u)(x,t)$, then (23) yields that 
   \begin{eqnarray}
   e^{-C_2\rho s}e(u)\leq C_4(1+\rho)e^{C_1s}\epsilon, \nonumber
   \end{eqnarray}
   $\forall (x,t)\in M\times [t_0,T]$. Take the supremum on the left, we have
   \begin{eqnarray}
   \frac{e^{-C_2\rho s}\rho}{1+\rho}\leq C_4e^{C_1s}\epsilon. \nonumber
   \end{eqnarray}
   Define the function $\phi$ by 
   \begin{eqnarray}
   \phi(x)=\frac{e^{-C_2sx}x}{1+x}. \nonumber
   \end{eqnarray}
   Compute that 
   \begin{eqnarray}
   \phi^\prime(x)&=&\frac{(-C_2se^{-C_2sx}x+e^{-C_2sx})(1+x)-e^{-C_2sx}x}{(1+x)^2} \nonumber\\
   &=&\frac{e^{-C_2sx}}{(1+x)^2}(-C_2sx^2-C_2sx+1), \nonumber
   \end{eqnarray}
   then $\phi$ attach its maximal value at maximal point $x_0=-\frac{1}{2}+\sqrt{\frac{1}{4}+\frac{1}{C_2s}}$ with value 
   $\phi(x_0)$. If $\epsilon$ is small enough such that
   \begin{eqnarray}
   C_4e^{C_1s}\epsilon < \phi(x_0), \nonumber
   \end{eqnarray}
   we can find that 
   \begin{eqnarray}
   \rho < -\frac{1}{2}+\sqrt{\frac{1}{4}+\frac{1}{C_2s}} \nonumber
   \end{eqnarray}
   or
   \begin{eqnarray}
   \rho > -\frac{1}{2}+\sqrt{\frac{1}{4}+\frac{1}{C_2s}}. \nonumber
   \end{eqnarray}
   Since for $T=t_0$
   \begin{eqnarray}
   \rho\leq 2D< -\frac{1}{2}+\sqrt{\frac{1}{4}+\frac{1}{C_2s}} \nonumber
   \end{eqnarray}
   and $\rho$ is a continuous function of $T$, we have $\rho < -\frac{1}{2}+\sqrt{\frac{1}{4}+\frac{1}{C_2s}}$. Finally we conclude that
   \begin{eqnarray}
   \rho < -\frac{1}{2}+\sqrt{\frac{1}{4}+\frac{1}{C_2s}} ,
   \end{eqnarray}
   for all $t$. Then there exits a sequence $t_i$, such that 
   \begin{eqnarray}
   u(\cdot,t_i)\rightarrow u_{\infty},    
   \end{eqnarray}
   where $u_{\infty}$ is a Lipschitz map. By the regularity of the
   solution of the subelliptic parabolic system (cf.\cite{MR3844509}), we only need  to show that $|u_{t_i}|\rightarrow 0$ as $t_i\rightarrow \infty,\ \forall x\in  M$.\par
   First note that, from (20) and (26)  we have 
   \begin{eqnarray}
   \int_{t_i-1}^{t_i}\int_M |u_r|^2(x,r)dVdr=E_b(t_i-1)-E_b(t_i)\rightarrow 0 
   \end{eqnarray}
   as $t_i\rightarrow \infty$. Since 
   \begin{eqnarray}
   (\triangle_b-\partial_t)|u_t|^2&=&2(u_{ti})^2+2u_{tii}u_t-2u_{tt}u_t  \nonumber\\
   &\geq& -2\kappa e_b(u)|u_t|^2\geq -2\kappa\rho|u_t|^2, \nonumber
   \end{eqnarray}
   where $\rho=\sup\limits_{M\times[0,\infty)}e(u)$ and $K^N\leq
   \kappa$. Then we find that
   \begin{eqnarray}
   (\triangle_b-\partial_t)e^{-2\kappa \rho t}|u_t|^2 \geq 0, \nonumber
   \end{eqnarray}
   hence applying Lemma 2.3, we conclude that 
   \begin{eqnarray}
   e^{-2\kappa \rho t_i}|u_t|^2(x,t_i)&\leq& C\int_{t_i-1}^{t_i}\int_Me^{-2\kappa \rho r}|u_r|^2(x,r)dVdr \nonumber \\
   &\leq&Ce^{-2\kappa \rho(t_i-1)}\int_{t_i-1}^{t_i}\int_M|u_r|^2(x,r)dVdr .\nonumber
   \end{eqnarray}
   Together with (27), this implies that
   \begin{eqnarray}
   |u_t|^2(x,t_i)&\leq&Ce^{2\kappa \rho}\int_{t_i-1}^{t_i}\int_M|u_r|^2(x,r)dVdr \nonumber\\
   &&\rightarrow 0 ,\nonumber
   \end{eqnarray}
   as $t_i\rightarrow \infty$. This completes the proof.\\
   {\qed}\\
   
   \section{Achnowledgments}
   
   \ \ \ The author would like to exprss his thanks to Professor Yuxin Dong and Professor Yibin Ren for valuable discussions and helpful suggestions.

\end{document}